\documentclass[12pt]{amsart}
\font\emailfont=cmtt10

\usepackage{amsmath,amsthm,amsfonts,amscd,flafter,epsf}

\title[{On Park's exotic smooth four-manifolds}] 
{On Park's exotic smooth four-manifolds}

\author[Peter Ozsv{\'a}th]{Peter Ozsv\'ath}
\address{Department of
Mathematics, Columbia University, New York 10027 \newline
\indent{\emailfont{petero@math.columbia.edu}}}
\thanks{PSO was supported by NSF grant number DMS 0234311}

\author[Zolt{\'a}n Szab{\'o}]{Zolt{\'a}n Szab{\'o}} 
\address{Department of
Mathematics, Princeton University, New Jersey 08544 \newline
\indent{\emailfont{szabo@math.princeton.edu}}}
\thanks{ZSz was supported by NSF grant number DMS 0107792}

\newcommand\mCP{\overline {\mathbb{CP}}}
\newcommand\PD{\mathrm{PD}}

\newtheorem{theorem}{Theorem}[section]

\newtheorem{lemma}[theorem]{Lemma}

\def\endproof{\relax\ifmmode\expandafter\endproofmath\else
  \unskip\nobreak\hfil\penalty50\hskip.75em\hbox{}\nobreak\hfil\bull
  {\parfillskip=0pt \finalhyphendemerits=0 \bigbreak}\fi}
\def\endproofmath$${\eqno\bull$$\bigbreak}
\def\bull{\vbox{\hrule\hbox{\vrule\kern3pt\vbox{\kern6pt}\kern3pt\vrule}\hrule}}

\newcommand{\Q}{\mathbb{Q}}
\newcommand{\R}{\mathbb{R}}

\newcommand{\Z}{\mathbb{Z}}

\newcommand{\CP}[1]{{\mathbb{CP}}^{#1}}

\newcommand{\cm}{\cdot}

\newcommand{\SW}{SW}

\newcommand{\SpinC}{{\mathrm{Spin}}^c}

\newcommand\sgn{\mathrm{sgn}}

\newcommand\abuts\Rightarrow

\begin{document}

\begin{abstract}  
	In a recent paper, Park constructs certain exotic
	simply-connected four-manifolds with small Euler
	characteristics. Our aim here is to prove that the
	four-manifolds in his constructions are minimal.
\end{abstract}

\maketitle
\section{Introduction}

Since the seminal works of Donaldson~\cite{DonaldsonChambers} and Freedman~\cite{Freedman}, it has
been known that closed, simply-connected four-manifolds can support
exotic smooth structures. In fact, for many homeomorphism classes,
gauge theory tools (Donaldson invariants and Seiberg-Witten
invariants) have been very successful at proving the existence of
infinitely many smooth structures, see for
example~\cite{DonaldsonPolynomials}, \cite{GompfMrowka},
\cite{FriedmanMorgan}, \cite{FSKnots}. However, exotic examples with
small Euler characteristics are much more difficult to find. For a
long time, the smallest known example was the Barlow
surface~\cite{Kotschick}, which has Euler characteristic $11$ and
which is homeomorphic, but not diffeomorphic, to $\CP{2}\# 8 \mCP^2$.
Recently, in a remarkable paper, Park~\cite{Park} constructs a symplectic
manifold $P$ with Euler characteristic $10$ using the rational
blow-down operation of Fintushel and Stern~\cite{QBD}, and proves that
it is homeomorphic, but not diffeomorphic to $\CP{2}\# 7\mCP^2$.

In this note, we compute the Seiberg-Witten invariants of $P$ and
prove the following:

\begin{theorem}
\label{thm:Park1}
Park's example $P$ does not contain any smoothly embedded two-spheres with
self-intersection number $-1$; equivalently, it is not the blow-up of
another smooth four-manifold.
\end{theorem}

In a similar manner, Park also constructs a symplectic four-manifold $Q$
which is  homeomorphic, but not diffeomorphic, to $\CP{2}\# 8\mCP^2$.
We prove here the following:

\begin{theorem}
\label{thm:Park2}
The manifold $Q$ contains no smoothly embedded two-sphere with
self-intersection number $-1$, and in particular $Q$ is not
diffeomorphic to $P\# \mCP^2$.
\end{theorem}

Note that $Q$ and the Barlow surface have the same Seiberg-Witten invariants,
and we do not know whether or not they are diffeomorphic.

\noindent{\bf Acknowledgements} The authors wish to thank Jongil Park, Jacob Rasmussen, and Andr{\'a}s Stipsicz
for interesting conversations during the course of this work.

\section{Seiberg-Witten theory}

We will deal in this paper with Seiberg-Witten theory for
four-manifolds $X$ with $b_2^+(X)=1$ (and $b_1(X)=0$). For the
reader's convenience, we recall the basic aspects of this theory, and
refer the reader to~\cite{KMthom}, \cite{Morgan} for more in-depth
discussions.

The Seiberg-Witten equations can be written down on any four-manifold
equipped with a $\SpinC$ structure and a Riemannian metric. We
identify here $\SpinC$ structures over $X$ with characteristic classes
for the intersection form of $X$, by taking the first Chern class of
the $\SpinC$ structure.  This induces a one-to-one correspondence in
the case where $H^2(X;\Z)$ has no two-torsion. Taking a suitable
signed count of solutions, one obtains a smooth invariant of $X$ when
$b_2^+(X)>1$. In the case where $b_2^+(X)=1$, the invariant depends on
the choice of the Riemannian metric through the cohomology class of
its induced self-dual two-form (compare
also~\cite{DonaldsonChambers}).

Formally, then, for a fixed two two-dimensional cohomology class $H\in
H^2(X;\R)$ with $H^2>0$ and characteristic vector $K\in H^2(X;\Z)$
with $K.H\neq 0$, the Seiberg-Witten invariant $\SW_{X,H}(K)$ is an
integer which is well-defined provided that $H.K\neq 0$.  This integer
vanishes whenever $$K^2< 2\chi(X)+3\sigma(X).$$ For fixed $H$, then,
the $H$-basic classes are those characteristic cohomology classes $K$
for which $\SW_{X,H}(K)\neq 0$.  The quantity
$K^2-2\chi(X)-3\sigma(X)$ is four times the formal dimension of the
moduli space of solutions to the Seiberg-Witten equations over $X$ in
the $\SpinC$ structure whose first Chern class is $K$. The
Seiberg-Witten invariant vanishes when this formal dimension is
negative; when it is positive, one cuts down the moduli space by a
suitable two-dimensional cohomology class to obtain an integer-valued
invariant.

More precisely, a Riemannian metric on $X$ induces a Seiberg-Witten
moduli space. The signed count of the solutions in this moduli space
depends only on the cohomology class of the induced self-dual two-form
$\omega_g$, which in the above case was denoted by $H$.  The
dependence on $H$ is captured by the wall-crossing
formula~\cite{KMthom}, \cite{LiLiu}: if $X$ is a four-manifold with
$b_1(X)=0$, and $H$ and $H'$ are two cohomology classes with
positive square and $H.H'>0$,
then $$\SW_{X,H}(K)=\SW_{X,H'}(K) +\left\{\begin{array}{ll} 0
&{\text{if $K.H$ and $K.H'$ have the same sign}}
\\
\pm 1 &{\text{otherwise.}}
\end{array}\right.
$$

It follows readily from the compactness result for the moduli space of
solutions to the Seiberg-Witten equations that for any $H$, there are
only finitely many $H$-basic classes.

It is interesting to note that the wall-crossing formula together with
the dimension formula (which states that $\SW_{X,H}(K)=0$ when
$K^2-2\chi(X)-3\sigma(X)<0$), ensures that if $X$ is a four-manifold
with $b_2^+(X)=1$ but $b_2(X)\leq 9$, there is only one chamber.

\subsection{Rational blow-downs}

In~\cite{QBD}, Fintushel and Stern introduce a useful
operation on smooth four-manifolds, and calculate how the
Seiberg-Witten invariants transform under this
operation. Specifically, let $C_p$ be the four-manifold which is a
regular neighborhood of a chain of two-spheres $\{S_0,...,S_p\}$ where
$S_0$ has self-intersection number $-4-p$, and $S_i$ has
self-intersection number $-2$ for all $i>0$.  The boundary of this
chain (the lens space $L((p+1)^2, p)$) also bounds a four-manifold
$B$ with $H^2(B;\Q)=0$.  If $X$ is a smooth, oriented four-manifold
with $b_2^+(X)>1$ which contains $C_p$, then we can trade $C_p$ for
the rational ball $B$ to obtain a new four-manifold $X'$. Clearly,
$H^2(X')$ is identified with the orthogonal complement to
$[S_i]_{i=0}^p$ in $H^2(X)$.

For each $\SpinC$ structure over $L((p+1)^2,p+1)$ which extends over
$B$, there is an extension (as a characteristic vector $K_0$) over
$C_p$ with the property that $K_0^2-p-1=0$.

Fintushel and Stern show that for any characteristic vector $K$ for
the intersection form of $X'$, $$\SW_{X'}(K)=\SW_{X}({\widetilde
K}),$$ where ${\widetilde K}$ is obtained from $K$, by extending over
the boundary by the corresponding characteristic vector $K_0$ as
above.

In the case where $b_2^+(X)=1$, the relation is expressed by choosing
a chamber for $X$ (and induced chamber for $X'$) whose metric form $H$
is orthogonal to each sphere in the configuration $C_p$.

\section{The four-manifold $P$}

We review Park's construction of $P$ briefly.  Start with a rational
elliptic surface with an ${\widetilde{E_6}}$ singularity (a
configuration of $-2$ spheres arranged in a star-like pattern, with a
central node and three legs of length two). There is a model of the
rational elliptic surface with the property that there are four nodal curves
in a complement of this singularity.  Blowing up the nodal curves, one
obtains four spheres of square $-4$. A section of the rational
elliptic surface meets all four of these spheres, and also one of the
leaves in the ${\widetilde{E_6}}$ singularity. Adding the section and
the four $-4$-spheres, one obtains a sphere $R_0$ with
self-intersection number $-9$ and then inside the ${\widetilde{E_6}}$
singularity, this can be extended to a chain of embedded spheres with
self-intersection $-2$ $\{R_i\}_{i=1}^5$. Park's example $P$ is
obtained by performing a rational blow-down, in the sense of Fintushel
and Stern~\cite{QBD}, on the chain of spheres $\{R_i\}_{i=0}^5$.
Since the spheres are all symplectic, a result of
Symington~\cite{Symington} guarantees that $P$ is symplectic.  

Theorem~\ref{thm:Park1} follows from the following refinement:

\begin{theorem}
Let $K$ denote the canonical class of $P$. Then, the Seiberg-Witten basic
classes of $P$ are $\{\pm K\}$.
\end{theorem}

It follows at once that $X$ is minimal. Specifically, if one could
write $X\cong Y \#\mCP^2$, then according to the blow-up
formula~\cite{FSthom}, the basic classes of $X$ come in pairs of the
form $K_0\pm E$ where $K_0$ runs over the basic classes of $Y$, and
$E$ denotes the exceptional curve in $\mCP^2$. But this is impossible
since $K^2=2$.

\begin{figure}
\mbox{\vbox{\epsfbox{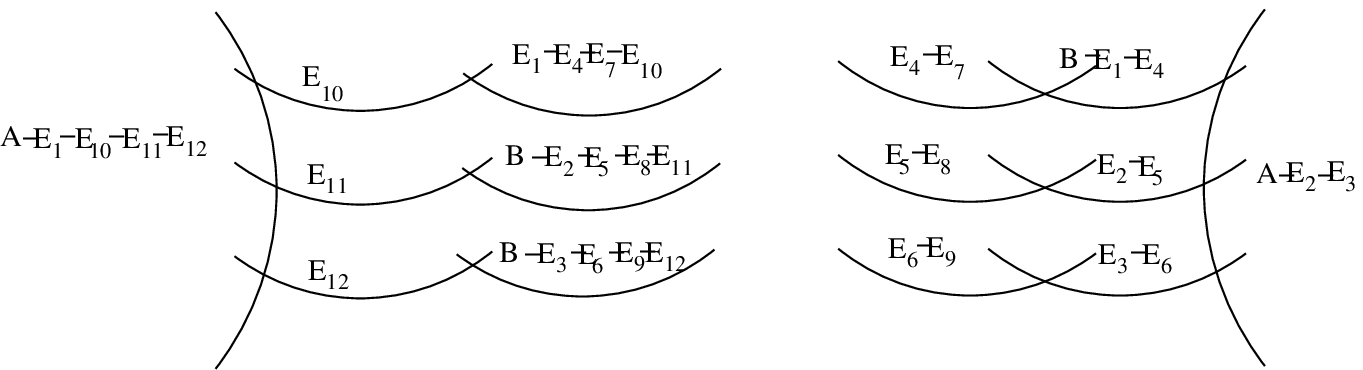}}}
\caption{\label{fig:ParkConfig}
We have illustrated here a basis of two-spheres for $\CP{2}\# 12\mCP^2$.}
\end{figure}

We find it convenient to describe the manifold $P$ in a concrete model.
Specifically, consider the four-manifold $X=S^2\times S^2 \# 12\mCP^2$,
with the basis of two-spheres $A$, $B$, $\{E_i\}_{i=1}^{12}$. Here,
$A$ and $B$ are supported in the $S^2\times S^2$ factor, so that
$A=\{a\}\times S^2$ and $B=S^2\times \{b\}$, while $E_i$ is the
``exceptional sphere'' (sphere of square $-1$) in the $i^{th}$ $\mCP^2$
summand. Alternatively, this manifold can be thought of as the blowup
of rational elliptic surface with an ${\widetilde{E_6}}$ singularity, and a
complementary singularity consisting of three $-1$-spheres arranged in
a triangular pattern, which is then blown up four times, to give a
tree-like configuration of spheres with a central sphere of of square
$-4$, and three legs consisting of a chain of a $-1$ sphere and
another $-4$ sphere. See Figure~\ref{fig:ParkConfig} for an illustration.

More precisely, consider the elliptic surface singularity which can be
described by three $-1$-framed unknots, each of which links the other
two in one point apiece. Denote the corresponding two-dimensional
homology classes by $A$, $B$, and $C$. It is well-known,
c.f.~\cite{HarerKasKirby} that this singularity can be perturbed into
four nodal curves. By blowing up the four double-points, we obtain
four disjoint $-1$-spheres. In fact, the homology class of the fiber
is represented by the homology class of the fiber $A+B+C$. Thus,
the four $-4$ spheres can written in the basis of homology as
$$\{A+B+C-2E_i\}_{i=1}^4,$$ where $E_i$ are the newly-introduced exceptional spheres.

Armed with this principle, the chain of spheres in $X$
which are to be
rationally blown down can be written homologically as:
\begin{eqnarray*}
R_0&=&10A+8B-6E_1-4E_2-4E_3-4E_4-4E_5-4E_6 \\ && -3E_7-4E_8-4E_9-2E_{10}-2E_{11}-2E_{12} \\
R_1&=&B-E_1-E_4 \\
R_2&=&A-E_2-E_3 \\
R_3&=&E_3-E_6 \\
R_4&=&E_6-E_9 \\
R_5&=&E_4-E_7
\end{eqnarray*}
Note that we are using here $E_7$ as our section, which is to be added
to the four $-4$-spheres coming from the complement of the ${\widetilde{E_6}}$
singularity.  The four exceptional spheres in the 
complementary singularity are represented by the spheres
$A-E_1$, $E_{10}$, $E_{11}$,
$E_{12}$.

Let $P$ denote the Park manifold obtained by rationally blowing down
the configuration $R_0,...,R_5$ in $X$. $\SpinC$ structures
over $P$ (labelled by characteristic vectors $K$)
correspond to characteristic vectors (labelled by characteristic vectors ${\widetilde K}$)
over $X$ whose evaluations on the configuration $\{R_i\}$ 
take one of the following
seven forms:
\begin{equation}
\label{eq:QBD}
\begin{array}{lrrrrr}
(7,& 0,& 0,& 0,& 0,& 0) \\
(-1,& 0,& -2,& 0,& 0,& 0) \\
(5,& 0,& 0,& 0,& 0,& -2) \\
(-3,& -2,& 0,& 0,& 0,& 0) \\
(3,& 0,& 0,& 0,& -2,& 0) \\
(-7,& 0,& 0,& 0,& 0,& 0) \\
(1,& 0,& 0,& -2,& 0,&  0) \\
\end{array}
\end{equation}
According to the rational blow-down formula~\cite{QBD},
$$\SW_{P}(K)=\SW_{X,H}({\widetilde K}),$$ where here $H\in H^2(X;\R)$
is any real two-dimensional cohomology class with $H^2>0$ and $H.H'>0$ and which is
orthogonal to all the $\{R_i\}$.  Moreover, according to the
wall-crossing formula, combined with the fact that
$S^2\times S^2$ has positive scalar curvature and hence
trivial invariants in a suitable chamber 
(c.f.~\cite{Witten}), it follows that 
$$\SW_{X,H}({\widetilde
K})=\left\{
\begin{array}{ll}
0 &{\text{if ${\widetilde K}^2+4<0$ or ${\widetilde K}.H$ and ${\widetilde K}.H'$ have the same sign}} \\
\pm 1 &{\text{otherwise,}}
\end{array}\right.
$$
where here $H'=\PD(A+B)$.
(The first condition for vanishing
is the dimension formula for the moduli space, while the second
condition comes from the wall-crossing formula.)

Explicitly, then, we see that the basic classes $K$ for $P$ are
precisely those for which the extension ${\widetilde K}$ (by one of
the vectors from the list in Equation~\eqref{eq:QBD}) satisfies:
${\widetilde K}^2+4\geq 0$ and also $\sgn({\widetilde K}.H)\neq
\sgn({\widetilde K}.H')$, where here $H$ is any (real) cohomology
class with $H^2>0$ and $H.H'>0$ and which is orthogonal to all the
$\{R_i\}_{i=0}^5$.  For example, we could use the vector $$H=(105, 92,
-67, -51, -41, -38, -36, -41, -38, -36, -41, -18, -18, -18)$$
(written here with respect to the basis Poincar\'e dual to
$\{A,B,E_1,...,E_{12}\}$).
In order to make this a finite computation, we proceed as follows.

Suppose that $Z$ is a smooth four-manifold with $b_2^+(Z)>1$, and we
have homology classes ${\mathbf C}=\{C_i\}_{i=1}^n$ with negative
self-intersection number $C_i\cm C_i = -p_i<0$. A cohomology class
$K\in H^2(X;\Q)$ is called ${\mathbf C}$-adjunctive if for each $i$
$\langle K,[C_i]\rangle$ is integral, and indeed the following two
conditions are satisfied:
\begin{eqnarray*}
|\langle K, [C_i] \rangle |&\leq& p_i  \\
\langle K,[C_i]\rangle &\equiv& p_i\pmod{2}.
\end{eqnarray*}
Clearly, the set of ${\mathbf C}$-adjunctive cohomology classes 
has size $\prod_{i=1}^n (p_i+1)$.

\begin{lemma}
	\label{lemma:Adjunctive}
	Let ${\mathbf S}=\{S_i\}_{i=1}^n$ be a collection of embedded
	spheres in $X$ whose homology classes are orthogonal to the
	the $\{R_i\}_{i=0}^5$. Let ${\mathbf C}=\{C_i\}_{i=1}^8$
	denote their induced homology classes in $H_2(P)$. If every
	${\mathbf C}$-adjunctive basic class for $P$ is
	zero-dimensional, then in fact every basic class for $P$ is
	${\mathbf C}$-adjunctive.
\end{lemma}

\begin{proof}
If $P$ has a basic class $L_0$ which is not ${\mathbf C}$-adjunctive, then
by the rational blow-down formula, 
$X$ has a basic class $L_1$
and a smoothly embedded sphere $S_i$ for which $|\langle L_1,[S_i]\rangle|> -S_i\cdot S_i$, where
we can use any metric whose period point $H'$ is perpendicular to the configuration
$\{R_i\}_{i=0}^5$. By fixing $H'$ to be also perpendicular to $S_i$, and using
the adjunction formula  for spheres of negative square~\cite{FSthom} we get
another basic class $L_2=L\pm 2\PD[S_i]$ of $X$. Applying the blowdown formula once more we get a
basic class $L_3$ of $P$ where the dimension of $L_3$ is bigger then the dimension
of $L_0$. Since $P$ has only finitely many basic classes this process has to stop, which means that
the final $L_{3k}$ class is ${\mathbf C}$-adjunctive. However it is also positive dimensional, thus
proving the lemma.
\end{proof}

Our next goal, then is to find a collection of embedded spheres
$\{S_i\}_{i=1}^{8}$ in $X$ which, together with the $\{R_i\}_{i=0}^5$
form a basis for
$H^2(X;\Q)$. To this end, we use the spheres:
\begin{eqnarray*}
S_1 &=& E_5-E_8 \\
S_2 &=& E_{12}-E_{10} \\
S_3 &=& E_{11}-E_{12} \\
S_4 &=& A-E_1-E_{11} \\
S_5 &=& A+B-E_1-E_2-E_5-E_8 \\
S_6 &=& -E_5+E_{10}+E_{11} \\
S_7 &=& 2 E_7+2E_4 -2A+E_{11} \\
S_8 &=& E_6+E_9+E_3-E_2 -2E_5.
\end{eqnarray*}
The spheres $\{S_i\}_{i=1}^5$ have square $-2$, while $S_6$, $S_7$, and $S_8$
have squares $-3$, $-9$, and $-8$ respectively.
It is easy to see that these classes are all orthogonal to the homology
classes generated by the spheres $\{R_i\}_{i=0}^5$.

It is easy to see, now, that there are $612360$ $\{S_i\}$-adjunctive
vectors in $H^2(X;\Q)$ with integral evaluations on each of the $S_i$,
and whose extension over the blow-down configuration is one of the
seven choices enumerated in Equation~\eqref{eq:QBD}. Of these, $12498$
correspond to characteristic cohomology classes in $H^2(X;\Z)$. Of
these, $8960$ have length $\geq -4$ (i.e. satisfying
$K^2-(2\chi+3\sigma)\geq 0$). Finally, only two of these have the
property that evaluation of $H$ and $H'$ have opposite sign. Indeed,
these classes are the canonical class $K$ and also $-K$. Since these
classes have dimension zero, it follows from
Lemma~\ref{lemma:Adjunctive} that these are the only two basic classes
for $P$.

\section{The four-manifold $Q$}

The manifold $Q$ is constructed as follows.  Start with a rational
surface with an ${\widetilde{E_6}}$ singularity as before, except 
now blow up only three of the nodes. In a manner similar to the
previous construction, one finds now a sphere of self-intersection
number $-7$ (gotten by resolving a section and the three $-4$
spheres). This is then completed by a chain of three $-2$ spheres in
the ${\widetilde{E_6}}$ singularity. Forming the rational blow-down,
one obtains a second manifold $Q$ which is homeomorphic to
$\CP{2}\#8\mCP^2$.

For $Q$, we prove the following:
\begin{theorem}
Let $K$ denote the canonical class of $Q$. Then, the Seiberg-Witten basic
classes of $Q$ are $\{\pm K\}$.
\end{theorem}

The second construction starts again with a rational surface.  For
this surface, we can take the previous one, only blow down the curve
$E_{12}$.

Again, we use the section $E_7$; now the three $-4$ spheres which are to be 
added are represented by $E_1-E_4-E_7-E_{10}$, $B-E_2-E_5-E_8-E_{11}$,
and $A-E_1-E_{10}-E_{11}$. Thus, our configuration which is to be rationally blown down consists of:
\begin{eqnarray*}
R_0&=&7A+6B-4E_1-3E_2-3E_3-3E_4-3E_5-3E_6\\
&&-2E_7-3E_8-3E_9-2E_{10}-2E_{11} \\
R_1&=&E_4-E_7 \\
R_2&=&B-E_1-E_4 \\
R_3&=&A-E_2-E_3.
\end{eqnarray*}

The vector 
$H=(229, 226, -143, -113, -113, -86, -87, -87, -86, -87, -87, -58, -58)$
has positive square, and is orthogonal to all the $\{R_i\}_{i=0}^3$.

A rational basis the cohomology of $(S^2\times S^2) \# 11\mCP^2$
is gotten by completing $R_0$, $R_1$, $R_2$, and $R_3$ with the following
set of spheres with negative square:
\begin{eqnarray*}
S_1&=&E_{10}-E_{11} \\
S_2&=&E_5-E_6 \\
S_3&=&E_8-E_9 \\
S_4&=&E_5-E_8 \\
S_5&=& E_2-E_3 \\
S_6 &=& A-E_1-E_{10} \\
S_7&=&A+B-E_1-E_2-E_5-E_8 \\
S_8&=&2A-2E_4-2E_7-E_{11} \\
S_9&=&2A+2B-E_1-E_2-E_3-E_4-E_7-2E_5-E_6-E_{10}.
\end{eqnarray*}

For this case, the $\SpinC$ structures over $L(25,4)$
which extend over the rational ball can be uniquely extended over
the configuration of spheres in one of the five possible ways:
\begin{equation}
\begin{array}{lrrr}
(5,& 0,& 0,& 0) \\
(-1,& -2,& 0,& 0) \\
(3, &0, &0, &-2) \\
(-5, &0, &0, &0) \\
(1,& 0,& -2,& 0).
\end{array}
\end{equation}

Again, there are $437400$ $\{S_i\}$-adjunctive vectors in $H^2$ with
rational coefficients, which have integral evaluations on each sphere
and which extend over the configuration of spheres $\{R_i\}_{i=0}^3$
as above. Of these, $17496$ correspond to (integral) characteristic
cohomology classes. Of these, 3754 have square $\geq -3$. Finally, of
these, exactly two ($K$ and $-K$) have the evaluations with opposite
sign against $H$ and $H'$, hence correspond to basic classes for
$Q$. Arguing as in Lemma~\ref{lemma:Adjunctive}, we see that these are
the only two basic classes for $Q$.

\bibliographystyle{plain}
\bibliography{biblio}
\end{document}